\documentclass[a4paper, 12pt]{article}
\usepackage{makeidx}
\usepackage{enumerate, theorem}
\usepackage{amsmath, amsfonts, amssymb}
\usepackage[height=22.5cm, width=15cm]{geometry}
\usepackage[all]{xy}
\usepackage{mathrsfs, yfonts}
\usepackage{graphicx}
\DeclareMathOperator{\C}{\mathbb{C}}
\newcommand{\parag}[1]{\paragraph{\sc{#1.}} }

\DeclareMathOperator{\id}{id}

\newtheorem{thm}{Theorem}[subsection]

\newtheorem{cor}[thm]{Corollary}
\newtheorem{prop}[thm]{Proposition}
\newtheorem{lemma}[thm]{Lemma}

\begin{document}

\title{Note on the semi-continuity of the algebraic dimension.}
\date{07/10/16}

 \author{Daniel Barlet\footnote{Institut Elie Cartan, G\'eom\`{e}trie,\newline
Universit\'e de Lorraine, CNRS UMR 7502   and  Institut Universitaire de France.}.}

 \maketitle
 
 \parag{Abstract} In this short Note we show  that the direct image sheaf $R^{1}\pi_{*}(\mathcal{O}_{\mathcal{X}})$ associated to an analytic family of compact complex manifolds $\pi : \mathcal{X} \to S$ parametrized by a reduced complex space $S$ is a locally free (coherent) sheaf of $\mathcal{O}_{S}-$modules. This result allows to improve a semi-continuity type result for the algebraic dimension of compact complex manifolds in an analytic family given in [B.15].
 
 \parag{AMS Classification  2010}  32G05-32A20-32J10. 
 
 \parag{Key words} Analytic family of compact complex manifolds. Algebraic dimension.
 
 \tableofcontents
 
 \newpage
 
 \section{Introduction} It is well known that for a compact complex manifold $X$ of the Fujiki-Varouchas class $\mathscr{C}$ (recall that Varouches [V.89] shows that this is simply the class of compact complex manifolds which admit a K\"ahler modification) the number
  $$h^{0, 1}(X) := \dim_{\C} H^{1}(X, \mathcal{O}_{X}) $$ 
  is a topological invariant (half of the first Betti number), so it is constant in an analytic deformation inside the class $\mathscr{C}$. In this short Note we prove that this number is invariant in any analytic deformation for any compact complex manifold. As an application, this allows us to improve a semi-continuity type result for the algebraic dimension of compact complex manifolds in an analytic family given in [B.15]. See the theorem \ref{second} and its corollaries in  section 3.
 
 \section{The result}
 
  \begin{thm}\label{utile 3}
 Let  \ $\pi : \mathcal{X} \to S$ \ be a holomorphic  family of compact complex connected manifolds of dimension \ $n$ \ parametrized by an irreducible complex space \ $S$. Then the coherent sheaf $R^{1}\pi_{*}(\mathcal{O}_{\mathcal{X}})$ on $S$ is a locally free sheaf. 
 \end{thm}
 
The proof will use several lemmata.
 
 \begin{lemma}\label{30}
 Let $X$ be a compact normal connected complex space of dimension $n$ and let $L$ be a holomorphic line bundle on $X$. Then if  $L$ and $L^{*}$ have a non trivial holomorphic  section,  the line bundle $L$ is holomorphically trivial.
 \end{lemma}
 
 \parag{proof} Let $\sigma$ and $\tau$ the non trivial holomorphic sections of $L$ and $L^{*}$. Then the function $x \mapsto < \sigma(x), \tau(x) >$ is holomorphic and not identically zero. So it is a non zero constant function and we see that $\sigma$ and $\tau$ cannot vanish. Now the map $L \to X \times \C$ given by $\xi \mapsto (\pi(\xi, \xi\big/ \sigma(\pi(\xi)))$ is holomorphic and linear on fibres with inverse the holomorphic map given by $(x, \lambda) \mapsto \lambda.\sigma(x)$. So $L$ is trivial.$\hfill \blacksquare$\\

 We shall use this lemma in order to get the fact that if a holomorphic line bundle on  a  compact complex connected manifold is not holomorphically trivial, then $L$ or $L^{*}$ has no non trivial holomorphic section.
  
 \begin{prop}\label{ouvert}
 Let $M$ be a reduced complex space and $(X_{s})_{s \in S}$ an analytic family of compact  $n-$cycles in $M$ parametrized by a closed  irreducible complex subset $S$ of $\mathcal{C}_{n}(M)$ the space of compact $n-$cycles in $M$. Assume that  each cycle in this family is reduced, normal and connected. Let $\mathcal{L}$ be a holomorphic  line bundle on $M$. Then the subset $\Sigma$  of points in $S$  such that the restriction $ \mathcal{L}_{\vert X_{s}}$ is not  holomorphically trivial is an open subset in $S$.
 \end{prop}
 
 Let $f : \mathcal{L} \to M$ be the projection. Then the direct image of compact $n-$cycles $f_{*} : \mathcal{C}_{n}(\mathcal{L}) \to \mathcal{C}_{n}(M)$ is a holomorphic map (see [B-M 1] chapter IV). We can restrict this map to the subset $Z \subset \mathcal{C}_{n}(\mathcal{L})$ defined by the condition that the cycles are connected and  that the direct image by $f$ of a cycle $C$ in $Z$ is a cycle $X_{s}$ for some $s \in S$. These two conditions are analytic and closed thanks to the theorem  IV 7.2.1 of [B-M 1] and to the assumption that $S$ is a closed analytic subset in $\mathcal{C}_{n}(M)$.\\
 Remark that, as we assume that $X_{s}$ is normal and connected, a compact connected $n-$cycle $C$ in $\mathcal{L}$ with direct image $X_{s}$ by the projection $f$ is a section of the line bundle $\mathcal{L}_{\vert X_{s}}$.\\
 Note that we have a closed embedding $ j : S \to Z$ which associates to $s \in S$ the reduced $n-$cycle in $Z$ equal to the zero section of the line bundle $\mathcal{L}_{\vert X_{s}}$.\\
 Now if  the line bundle $\mathcal{L}_{\vert X_{s}}$ has a non trivial holomorphic section  the cycle $j(s)$ can move in $Z \cap f_{*}^{-1}(j(s))$ by homotheties in an analytic $1-$dimensional family containing $j(s)$. So we have
 $$ \dim_{j(s)} \big(Z \cap f_{*}^{-1}(j(s))\big) \geq 1.$$
 But the subset $W$ of points $w$ in $Z$ such the inequality $\dim_{w} [Z \cap f_{*}^{-1}(f_{*}(w))] \geq 1$ is a closed analytic subset in $Z$. So the subset $\Sigma_{0} : = j^{-1}(W)$ is a closed analytic subset in $S$. Then the complement of $\Sigma_{0}$ is an open set in $S$. So if $L_{\vert X_{0}}$ has no non trivial holomorphic section, for $s$ in this open set, $L_{\vert X_{s}}$ is not holomorphically trivial. If $L^{*}_{\vert X_{0}}$ has no non trivial holomorphic section we obtain in the same way an open set around $0$ such that, for any $s$ in it,  $L_{\vert X_{s}}$ is not holomorphically trivial. The case when $L$ and $L^{*}$ have both a non trivial holomorphic section is excluded by the lemma \ref{30}.$\hfill \blacksquare$

 \begin{lemma}\label{manquant 2}
 Let $\pi : \mathcal{X} \to \Delta$ a proper holomorphic submersion of a complex manifold $\mathcal{X}$ onto an open  disc $\Delta$ with center $0$  in $\C$, with $n-$dimensional connected fibres. Let $L$ be a line bundle on $\mathcal{X}$ and assume that $L$ is holomorphically trivial on each $X_{s}, \forall s \in \Delta$. Then $L$ is trivial on $\mathcal{X}$.
  \end{lemma}
  
  \parag{proof} Consider the following data : an open disc $\Delta_{1}$ in $\Delta$, an open set $\mathcal{U}$ in $\pi^{-1}(\Delta_{1})$, a holomorphic  trivialization $t : L_{\vert \mathcal{U}} \to \mathcal{U}\times \C$ and a holomorphic section $ \gamma : \Delta_{1} \to \mathcal{U}$ of $\pi$. Of course, choosing first a local trivialization of $L$ on an open set in $\mathcal{X}$ we can find such data with any point in $\Delta$ as the center of the (small) disc $\Delta_{1}$.\\
  Let $Z \subset \mathcal{C}_{n}(L)$ the analytic subset of connected compact $n-$cycles $C$  in $L$ such the direct image cycle $f_{*}(C)$ of $C$ by the projection $f : L \to \mathcal{X}$ is one of the fibres of $\pi$. So we have a holomorphic map $g : Z \to \Delta$ defined by $ f_{*}(C) = X_{g(C)}$. Denote now by $Z_{1}$ the subset in $Z$ of cycles $C \in g^{-1}(\Delta_{1})$ such that  the cycle $t_{*}(C\cap f^{-1}(\mathcal{U}))$  contains the point $(\gamma(g(C)), 1) \in \mathcal{U} \times \C$. we want to prove the following assertions :
  \begin{enumerate}[1)]
  \item The subset $Z_{1}$ is a closed analytic subset of  the open set $g^{-1}(\Delta_{1}) \subset Z$.
  \item The projection on $L_{\vert \pi^{-1}(\Delta_{1})}$ of the  graph $\Gamma_{1} \subset Z_{1}\times L_{\vert \pi^{-1}(\Delta_{1})} $ of the analytic family of compact $n-$cycles in $L$ parametrized by $Z_{1}$ is a closed embedding of a  complex sub-manifold in $L_{\vert \pi^{-1}(\Delta_{1})}$ which is disjoint of the zero section and gives a holomorphic section of  $L_{\vert \pi^{-1}(\Delta_{1})}$.
  \end{enumerate}
  As a consequence, we shall obtain that $L_{\vert \pi^{-1}(\Delta_{1})}$ is trivial on $\pi^{-1}(\Delta_{1})$. And, as this is true for any given point $s_{1}$ in $\Delta$ and a small enough open disc $\Delta_{1}$ with center $s_{1}$, the conclusion will follow because $H^{1}(\Delta, \mathcal{O}_{\Delta}^{*}) = \{1\}$.\\
  
  Let us prove the assertion 1).
  As the condition for $C \in g^{-1}(\Delta_{1})$ to be in $Z_{1}$ is given by the fact that the point $(C, (\gamma(g(C)), 1))$ is in the image of the graph $\Gamma_{1}\cap (Z_{1}\times L_{\vert\mathcal{U}})$ by the proper embedding $\id_{Z_{1}}\times\, t$, this is clearly a closed analytic condition as $g, \gamma$ and $t$ are holomorphic.\\
  
  To prove the assertion 2), remark first that each $C \in Z_{1}$ is the image of a holomorphic section of $L_{\vert X_{g(C)}}$ which does not vanishes at the point $\gamma(g(C))$. As $L_{\vert X_{g(C)}}$ is trivial, this section never vanishes on $X_{g(C)}$. Remark also that $g$ is injective in $Z_{1}$ because if $g(C) = g(C') := s$ then $C$ and $C'$ in $Z_{1}$ are the images of two holomorphic sections of the trivial line bundle $L_{\vert X_{s}}$ and take the same value at the point $\gamma(s) $. So $C = C'$ and $g : Z_{1}\to \Delta_{1}$ is an isomorphism. So the analytic family of compact cycle $(C)_{C \in Z_{1}}$ gives exactly one holomorphic never vanishing  section of $L_{\vert X_{s}}$ for each $s \in \Delta_{1}$. This is enough to prove our second assertion as the graph of this analytic family is a closed analytic subset in $L_{\pi^{-1}(\Delta_{1})}$ disjoint from the zero section and which is one to one on $\pi^{-1}(\Delta_{1})$ by the projection of $L$ on $\mathcal{X}$.$\hfill \blacksquare$\\

 \begin{lemma}\label{manquant}
  Let $\pi : \mathcal{X} \to \Delta$ a proper holomorphic submersion of a complex manifold $\mathcal{X}$ onto an open  disc $\Delta$ with center $0$  in $\C$, with $n-$dimensional connected fibres. Consider the injection of sheaves on $\Delta$
$$  j : R^{1}\pi_{*}\mathbb{Z} \to R^{1}\pi_{*}\mathcal{O}_{\mathcal{X}}.$$
The following properties are equivalent:
\begin{enumerate}[1)]
\item  Any section $\sigma$ with support $\{0\}$ of the sheaf $R^{1}\pi_{*}\mathcal{O}_{\mathcal{X}}$ vanishes.
\item Any section $\sigma$ of the sheaf $R^{1}\pi_{*}\mathcal{O}_{\mathcal{X}}$ such its restriction to $\Delta^{*}$ is   in the image of the map
 $j : H^{0}(\Delta^{*}, R^{1}\pi_{*}\mathbb{Z}) \to H^{0}(\Delta^{*} ,R^{1}\pi_{*}\mathcal{O}_{\mathcal{X}}) $ is also in the image of  the map
 $j : H^{0}(\Delta, R^{1}\pi_{*}\mathbb{Z}) \to H^{0}(\Delta, R^{1}\pi_{*}\mathcal{O}_{\mathcal{X}}) $.
\item Any topologically trivial line bundle on $\mathcal{X}$ which induces on $X_{0}$  a line bundle which is  holomorphically trivial on each $X_{s}$ for any $s \not= 0$ near-by enough $0$  induces a line bundle which is   holomorphically trivial on $X_{0}$.
\end{enumerate}
\end{lemma}

\parag{proof} $1) => 2)$. Take any section $\sigma$ of the sheaf $R^{1}\pi_{*}\mathcal{O}_{\mathcal{X}}$ such its restriction to $\Delta^{*}$ is in the image of  $j : H^{0}(\Delta^{*}, R^{1}\pi_{*}\mathbb{Z}) \to H^{0}(\Delta^{*} ,R^{1}\pi_{*}\mathcal{O}_{\mathcal{X}}) $. As $R^{1}\pi_{*}\mathbb{Z}$ is a constant sheaf on $\Delta$ we have $H^{0}(\Delta^{*}, R^{1}\pi_{*}\mathbb{Z}) = H^{0}(\Delta, R^{1}\pi_{*}\mathbb{Z})$. So there exists $\tau \in  H^{0}(\Delta, R^{1}\pi_{*}\mathbb{Z})$ such that $\sigma - j(\tau)$ vanishes on $\Delta^{*} $. Then by 1) we have $\sigma = j(\tau)$.\\

$2) => 3)$. As $\Delta$ is Stein and contractible and we know that  for each $i \geq 0$ the sheaves $R^{i}\pi_{*}\mathcal{O}_{\mathcal{X}} $ are coherent and the sheaves $R^{i}\pi_{*}(\mathbb{Z})$ are constant sheaves, the Leray spectral sequence gives natural  isomorphisms $H^{i}(\mathcal{X}, \mathcal{O}_{\mathcal{X}}) \simeq H^{0}(\Delta, R^{i}\pi_{*}\mathcal{O}_{\mathcal{X}} )$ and $H^{i}(\mathcal{X}, \mathbb{Z}) \simeq H^{0}(\Delta,  R^{i}\pi_{*}(\mathbb{Z}))$ for each $i \geq 0$. Then we have:

$$ Coker j := H^{0}(\Delta, R^{1}\pi_{*}\mathcal{O}_{\mathcal{X}})\big/j( H^{0}(\Delta, R^{1}\pi_{*}\mathbb{Z})) \simeq H^{1}(\mathcal{X}, \mathcal{O}_{\mathcal{X}})\big/H^{1}(\mathcal{X}, \mathbb{Z}) $$
which classifies the holomorphic line bundles on $\mathcal{X}$ which are topologically trivial, up to isomorphism. So the isomorphism class of a given  topologically trivial line bundle $L$  is defined by  the image in  $Coker j$ of some $\sigma \in H^{0}(\Delta, R^{1}\pi_{*}\mathcal{O}_{\mathcal{X}})$.\\
 Now take a line bundle $L$ on $\mathcal{X}$ which is topologically trivial. Assume that $L$ is holomorphically trivial on each $X_{s}$ for $s \in \Delta^{*}$.  So, thanks to the lemma \ref{manquant 2}, this  implies that the  section  $\sigma$  corresponding to the isomorphism class of $L$ is such that $ \sigma_{\vert \Delta^{*}}$ is in $j(H^{0}(\Delta^{*}, R^{1}\pi_{*}\mathbb{Z}))$. So by $2)$  we obtain that $\sigma$ gives $0$ in $Coker j$ and then the line bundle $L$ is holomorphically trivial. So the restriction to $X_{0}$ is holomorphically trivial and $3)$ is proved.\\

$3) => 1)$. If $\sigma \in H^{0}(\Delta, R^{1}\pi_{*}\mathcal{O}_{\mathcal{X}}) $ vanishes on $\Delta^{*}$ this implies that the corresponding line bundle on $\mathcal{X}$ is trivial on each $X_{s}, \forall s \in \Delta^{*}$. If $L_{X_{0}}$ is also trivial, then the lemma \ref{manquant 2} implies that $L$ is trivial on $\mathcal{X}$. So there exists some $\tau \in H^{0}(\Delta, R^{1}\pi_{*} \mathbb{Z})$ such  that $\sigma = j(\tau)$. But as $j$ is injective and as $R^{1}\pi_{*}\mathbb{Z}$ is a  constant sheaf, we have $\tau = 0$ and then $\sigma = 0$. So if $\sigma \not= 0$ the restriction $L_{\vert X_{0}}$ cannot be holomorphically trivial and then 3) gives a contradiction.$ \hfill \blacksquare$\\
 
  \begin{cor}\label{utile 1}
 Let $\pi : \mathcal{X} \to \Delta$ a proper holomorphic submersion of a complex manifold $\mathcal{X}$ onto an open  disc $\Delta$ with center $0$  in $\C$, with $n-$dimensional connected fibres. Then the coherent sheaf $R^{1}\pi_{*}(\mathcal{O}_{\mathcal{X}})$ is locally free.
 \end{cor}
 
 \parag{proof} It is enough to prove that the coherent sheaf $R^{1}\pi_{*}\mathcal{O}_{\mathcal{X}}$ has no torsion so that the property 1) in the previous lemma is satisfied. But he property 3) of the previous lemma is given by the proposition \ref{ouvert}. $\hfill \blacksquare$\\
  
  \parag{proof of the theorem \ref{utile 3}} It is enough to prove that this sheaf is $S-$flat. But the classical ``curve test'' for flatness\footnote{A geometric way to get this is to consider the linear space associated to this coherent sheaf : then on any curve it has  constant rank by corollary \ref{utile 1}; so it is a vector bundle and the sheaf is locally free.} is clearly satisfied thanks to the corollary \ref{utile 1}.$\hfill \blacksquare$\\

 \section{Application}
 
 As an immediate consequence of the theorem \ref{utile 3} we can suppress the hypothesis on the continuity of the $h^{0, 1}(s)$ in the theorem 1.0.3  of [B.15] and obtain the following semi-continuity result for the algebraic dimension.

 \begin{thm}\label{second}
Let  \ $\pi : \mathcal{X} \to S$ \ be a holomorphic  family of compact complex connected manifolds of dimension \ $n$ \ parametrized by an irreducible complex space \ $S$. Assume that there exists a dense Zariski open set \ $S'$ \ in \ $S$ \ such that for each \ $s$ \ in \ $S'$ \ the manifold \ $X_{s}$ \ satisfies the \ $\partial\bar \partial-$lemma\footnote{See for instance [Va.86].} and such that there exists a (smooth) relative sG-form for the family \ $\pi_{\vert S'} :\mathcal{X}_{\vert S'} \to S' $.\\
Then if \ $a : = \inf_{s \in S'}[a(X_{s})]$ \ we have \ $a(X_{s})\geq a$ \ for each \ $s \in S$.$\hfill \blacksquare$
\end{thm}

 \parag{Remark} A simpler statement (see remark 3 following the theorem 1.0.3 in [B.15]) which is a special case of the previous one, is obtained by assuming that the restriction of $\pi$ to $\pi^{-1}(S')$ is a weakly k\"ahler morphism in the sense of F. Campana (see for instance [C.81]); this implies the $\partial\bar \partial-$lemma assumption and the existence of a smooth relative sG-form for the restriction of $\pi$  over $S'$.$\hfill \square$\\
 
 As it is not so easy to show that a proper map is weakly K\"ahler (and we need less : each fibre in $S'$ has a sG-form and satisfies the $\partial\bar \partial-$lemma is enough)  let me recall the following results from [B.15]
 
 \begin{lemma}\label{Stab.}
Let \ $\pi : \mathcal{X} \to S$ \ be a proper holomorphic family of compact connected complex manifolds of dimension \ $n$ \ parametrized by an irreducible complex space \ $S$. Assume that for a point \ $s_{0} \in S$, the manifold  $X_{s_{0}} : = \pi^{-1}(s_{0})$ \ has a sG-form \ $\omega_{0}$. Then we can find a small open neighbourhood \ $S' $ \ of \ $s_{0}$ \ in \ $S$ \ and a relative sG-form \ $\omega$ \ on \ $\pi^{-1}(S')$ \ inducing \ $\omega_{0}$ \ on \ $X_{s_{0}}$.
\end{lemma}

 \begin{thm}\label{first}
Let  \ $\pi : \mathcal{X} \to S$ \ be a holomorphic  family of compact complex connected manifolds of dimension \ $n$ \ parametrized by an irreducible complex space \ $S$. Let  \ $s_{0}$ \ in \ $S$ \ such that  the manifold \ $X_{s_{0}}$ \ admits a (smooth) sG-form. Then there exists an open neighbourhood  \ $S_{0}$ \ of \ $s_{0}$,  a countable union \ $\Sigma$ \ of  closed irreducible analytic subsets in \ $S_{0}$ \ with no interior point  and a non negative integer \ $a $ \ such that  
\begin{enumerate}[(i)]
\item For any \ $s \in S_{0}$ \ we have \ $a(X_{s} )\geq a$.
\item For any \ $s \in S_{0}\setminus \Sigma $ \ we have \ $a(X_{s}) = a$.
\end{enumerate}
\end{thm}

Then the following corollaries are immediate from the theorem \ref{second} and \ref{first}.\\

\begin{cor}\label{1}
Let  \ $\pi : \mathcal{X} \to S$ \ be a holomorphic  family of compact complex connected manifolds of dimension \ $n$ \ parametrized by an irreducible complex space \ $S$. Assume that there exists a dense Zariski open set \ $S'$ \ in \ $S$ \ such that for each \ $s$ \ in \ $S'$ \ the manifold \ $X_{s}$ is K\"ahler. Then if \ $a : = \inf_{s \in S'}[a(X_{s})]$ \ we have \ $a(X_{s})\geq a$ \ for each \ $s \in S$.
\end{cor}

\bigskip

\begin{cor}\label{2}
Let  \ $\pi : \mathcal{X} \to S$ \ be a holomorphic  family of compact complex connected manifolds of dimension \ $n$ \ parametrized by an irreducible complex space \ $S$. Assume that there exists a dense Zariski open set \ $S'$ \ in \ $S$ \ such that for each \ $s$ \ in \ $S'$ \ the manifold \ $X_{s}$ is projective. Then for each $s \in S$ the manifold $X_{s}$ is Moishezon.
\end{cor}

We conclude by noticing that there exists an analytic family of smooth complex compact surfaces of the class VII (not K\"ahler) parametrized by a disc $\Delta$ such that the central fibre has algebraic dimension $0$ and all other fibres have algebraic dimension 1. See [F-P.09].\\
This shows that in our theorem \ref{second} some K\"ahler type assumption on the general fibre $X_{s}$ cannot be avoided in order that the ``general'' algebraic dimension gives a lower bound for the algebraic dimensions of all fibres.

   \bigskip

 \section{References}
 
 \begin{itemize}
 \item{[B.15]} Barlet, D. {\it Two semi-continuity results for the Algebraic Dimension of Compact Complex Manifolds} J. Math. Sci. Univ. Tokyo 22 (2015), pp.1-16.
 \item{[B-M 1]} Barlet, D. et Magn\`{u}sson, J. {\it Cycles analytiques complexes I, Th\'eor\`{e}mes de pr\'eparattion des cycles}, Cours Sp\'ecialis\'es 22, Soci\'et\'e Math\'ematique de France, Paris 2014.
 \item{[C.81]} Campana, F. {\it R\'{e}duction alg\'{e}brique d'un morphisme faiblement K{\"a}hl\'{e}rien propre et applications}, Math. Ann. 256 (1981), no. 2, pp.157-189.
 \item{[F-P.09]} Fujiki, A. and Pontecorvo, M. {\it Non-upper continuity of algebraic dimension for families of compact complex manifolds} \\
 arXiv: 0903.4232v2 [math. AG] .
 \item{[V.86]} Varouchas, J. {\it Propri\'{e}t\'{e}s cohomologiques d'une classe de vari\'{e}t\'{e}s analytiques complexes compactes}, S\'{e}minaire d'analyse P. Lelong-P. Dolbeault-H. Skoda,  1983/1984, pp.233-243, Lecture Notes in Math., 1198, Springer, Berlin, 1986.
 \item{[Va.89]} Varouchas, J. {\it K{\"a}hler spaces and proper open morphisms}, Math. Ann. 283 (1989), no. 1, pp.13-52.
  \end{itemize}

 \end{document}